\documentclass[a4paper,11pt]{amsart}
\usepackage{graphicx}
\usepackage{mathptmx}
\usepackage{amsmath}
\usepackage{amssymb}
\usepackage{enumitem}
\usepackage{xcolor}
\usepackage{xparse}
\NewDocumentCommand{\eulerian}{omm}
 {%
  \genfrac<>{0pt}{}{#2}{#3}%
  \IfValueT{#1}{_{\!#1}}%
 }

 \newmuskip\pFqmuskip

\newcommand*\pFq[6][8]{%
  \begingroup 
  \pFqmuskip=#1mu\relax
  \mathchardef\normalcomma=\mathcode`,
  \mathcode`\,=\string"8000
  \begingroup\lccode`\~=`\,
  \lowercase{\endgroup\let~}\pFqcomma
  {}_{#2}F_{#3}{\left(\genfrac..{0pt}{}{#4}{#5}\bigg|#6\right)}%
  \endgroup
}
\newcommand{\pFqcomma}{{\normalcomma}\mskip\pFqmuskip}

\newtheorem{theorem}{Theorem}
\newtheorem{lemma}[theorem]{Lemma}

\begin{document}

\title[dimorphic properties of Bernoulli random variable]{dimorphic properties of Bernoulli random variable}

\author{Taekyun  Kim}
\address{Department of Mathematics, Kwangwoon University, Seoul 139-701, Republic of Korea}
\email{tkkim@kw.ac.kr}

\author{DAE SAN KIM}
\address{Department of Mathematics, Sogang University, Seoul 121-742, Republic of Korea}
\email{dskim@sogang.ac.kr}

\author{Hyunseok  Lee}
\address{Department of Mathematics, Kwangwoon University, Seoul 139-701, Republic of Korea}
\email{luciasconstant@kw.ac.kr}

\author{Seong-Ho Park}
\address{Department of Mathematics, Kwangwoon University, Seoul 139-701, Republic of Korea}
\email{abcd2938471@kw.ac.kr}


\maketitle

\begin{abstract}
The aim of this paper is to study a dimorphic property associated with two different sums of identically independent Bernoulli random variables having two different families of probability mass functions.
In addition, we give two expressions on sums of products of degenerate Stirling numbers of the second kind and Stirling numbers of the first kind connected with those two different sums of identically independent Bernoulli random variables.
\end{abstract}

\section{Introduction}
It is well known that Bernoulli random variable is the discrete random variable which takes the value 1 with probability $p$ and the value 0 with probability $1-p$, where $0\le p\le1$.
In this paper, we study a dimorphic property (see Theorem 2) associated with two different sums of identically independent Bernoulli random variables having two different families of probability mass functions.\par
Further, we give two expressions on sums of products of degenerate Stirling numbers of the second kind and Stirling numbers of the first kind in connection with those two different sums of identically independent Bernoulli random variables. In fact, one is expressed in terms of the expectation of a random variable associated with one sum of identically independent Bernoulli random variables (see Theorem 3) and the other in terms of an integral involving the other sum of identically independent Bernoulli random variables (see Theorem 4). In the rest of this section, we recall some facts that are needed throughout this paper. \par
For any $\lambda\in\mathbb{R}$, the degenerate exponential function is defined as 
\begin{equation}
e_{\lambda}^{x}(t)=\sum_{n=0}^{\infty}(x)_{n,\lambda}\frac{t^{n}}{n!},\quad (\mathrm{see}\ [7,8,9,10,11,12]),\label{1}	
\end{equation}
where $(x)_{0,\lambda}=1,\ (x)_{n,\lambda}=x(x-\lambda)\big(x-(n-1)\lambda\big),\ (n\ge 1)$. \\
When $x=1$, for simplicity we write $e_{\lambda}(t)=e_{\lambda}^{1}(t)$. Note that $\displaystyle \lim_{\lambda\rightarrow 0}e_{\lambda}^{x}(t)=\sum_{n=0}^{\infty}\frac{x^{n}}{n!}t^{n}=e^{xt}\displaystyle$. \\ \par
The Stirling numbers of the first kind are defined by 
\begin{equation}
(x)_{n}=\sum_{l=0}^{n}S_{1}(n,l)x^{l},\quad(n\ge 0),\quad(\mathrm{see}\ [1,3,4,14]), \label{2}
\end{equation}
where $(x)_{0}=1,\ (x)_{n}=x(x-1)\cdots(x-n+1),\ (n\ge 1)$. 
As the inversion formula of \eqref{2}, the Stirling numbers of the second kind are defined by 
\begin{equation}
x^{n}=\sum_{l=0}^{n}S_{2}(n,l)(x)_{l},\quad(n\ge 0),\quad(\mathrm{see}\ [14,15,16,17,18]).\label{3}
\end{equation} \par
Moreover, in [8] the degenerate Stirling numbers of the second kind are defined as 
\begin{equation}
(x)_{n,\lambda}=\sum_{l=0}^{n}S_{2,\lambda}(n,l)(x)_{l},\quad(n\ge 0).\label{4}	
\end{equation}
Note that $\displaystyle \lim_{\lambda\rightarrow 0}S_{2,\lambda}(n,l)=S_{2}(n,l),\ (n,l\ge 0)\displaystyle$. From \eqref{4}, we note that 
 \begin{equation}
 \frac{1}{k!}\big(e_{\lambda}(t)-1\big)^{k}=\sum_{n=k}^{\infty}S_{2,\lambda}(n,k)\frac{t^{n}}{n!},\quad(k\ge 0),\quad(\mathrm{see}\ [6,8,9]).\label{5}
 \end{equation} \par
Let $X$ be a discrete random variable with probability mass function $p(j)=P\{X=j\},\ (j=1,2,\dots)$. Then the expectation of $X$ is defined by 
\begin{equation}
E[X^{n}]=\sum_{i=1}^{n}i^{n}p(i),\quad(n\in\mathbb{N}),\quad(\mathrm{see}\ [1,2,17]). \label{7}
\end{equation}
It is known that that the variance of $X$ is given by 
\begin{equation}
\sigma^{2}=\mathrm{Var}(X)=E[X^{2}]-\big(E[X]\big)^{2},\quad (\mathrm{see}\ [17]). \label{8}
\end{equation} \par
Let $(X_{j})_{1\le j\le n}$ be identically independent Bernoulli random variables such that $X_{j}$ has the probability of success $\frac{1}{j},\ (j=1,2,\dots,n)$. That is, 
\begin{equation}
X_{j}=\left\{\begin{array}{cc}
	1, & \textrm{if success,} \\
	0, & \textrm{otherwise,}
\end{array}	\right.\label{6}
\end{equation}
with $P\{X_{j}=1\}=\frac{1}{j},\ P\{X_{j}=0\}=1-P\{X_{j}=1\}$, where $j=1,2,3,\dots,n$, (see [2,5,13,17,18]). \par 
Let us assume that the random variable $Y_{n}$ is defined by 
\begin{equation}
Y_{n}=X_{1}+X_{2}+\cdots+X_{n}=\sum_{i=1}^{n}X_{i}. \label{9}	
\end{equation}
From \eqref{7}, \eqref{8} and \eqref{9}, we note that 
\begin{equation}
\mu_{n}=E[Y_{n}]=\sum_{j=1}^{n}E[X_{j}]=\sum_{j=1}^{n}\frac{1}{j}, \label{10}	
\end{equation}
and 
\begin{equation}
\sigma_{n}^{2}=E[Y_{n}^{2}]-\big(E[Y_{n}]\big)^{2}=\sum_{j=1}^{n}\frac{1}{j}\bigg(1-\frac{1}{j}\bigg).\label{11}	
\end{equation}

\section{Dimorphic properties of Bernoulli random variables} 
In this section, we show the dimorphic property (see Theorem 2) associated with the two different sums $Y_{n}=\sum_{j=1}^{n}X_{j}$ (see \eqref{9}) and $Z_{n,\lambda}(\alpha)=\sum_{j=1}^{n}X_{j,\lambda}(\alpha)$ (see \eqref{13}) of the identically independent Bernoulli random variables $X_{j}$, with the probability of success $\frac{1}{j}$, and those $X_{j,\lambda}(\alpha)$, with the probability of success $\frac{\alpha}{\lambda j+\alpha}$.\par
For $\alpha>0\ (\alpha\in\mathbb{R})$, and $\lambda\in(0,1)$, let $\big(X_{j,\lambda}(\alpha)\big)_{1\le j\le n}$ be identically independent Bernoulli random variables such that $X_{j,\lambda}(\alpha)$ has the probability of success $\frac{\alpha}{\lambda j+\alpha},\ (j=1,2,\dots,n)$. \\
That is,
\begin{equation}
X_{j,\lambda}(\alpha)=\left\{\begin{array}{cc}
	1, & \textrm{if success,} \\
	0 & \textrm{otherwise,}
\end{array}\right.	\label{12}
\end{equation}
with $P\{X_{j,\lambda}(\alpha)=1\}=\frac{\alpha}{\alpha+\lambda j},\ P\{X_{j,\lambda}(\alpha)=0\}=1-P\{X_{j,\lambda}(\alpha)=1\}$.\\
Let us assume that $Z_{n,\lambda}(\alpha)$ is defined by 
\begin{equation}
Z_{n,\lambda}(\alpha)=X_{1,\lambda}(\alpha)+X_{2,\lambda}(\alpha)+\cdots+X_{n,\lambda}(\alpha)=\sum_{i=1}^{n}X_{i,\lambda}(\alpha). \label{13}	
\end{equation}
From \eqref{7}, \eqref{8} and \eqref{13}, we have 
\begin{equation}
\mu_{n,\lambda}(\alpha)=E\big[Z_{n,\lambda}(\alpha)\big]=\sum_{j=1}^{n}\frac{\alpha}{\alpha+\lambda j},\label{14}	
\end{equation}
and 
\begin{align}
&\sigma_{n}^{2}(\alpha)\ =\ \mathrm{Var}\big(Z_{n,\lambda}(\alpha)\big)\ =\ E\big[Z_{n,\lambda}^{2}(\alpha)\big]-\Big(E\big[Z_{n,\lambda}(\alpha)\big]\Big)^{2}\label{15}\\
&=\ \sum_{j=1}^{n}E\big[X_{j,\lambda}^{2}(\alpha)\big]+2\sum_{i<j}E\big[X_{j,\lambda}(\alpha)\big]E\big[X_{i,\lambda}(\alpha)\big] \nonumber \\
&\quad -\bigg\{\sum_{j=1}^{n}\bigg(\frac{\alpha}{\alpha+\lambda j}\bigg)^{2}+2\sum_{i<j}\bigg(\frac{\alpha}{\alpha+\lambda j}\bigg)\bigg(\frac{\alpha}{\alpha+\lambda i}\bigg)\bigg\}\nonumber \\
&=\ \sum_{j=1}^{n}\frac{\alpha}{\alpha+\lambda j}+2\sum_{i<j}\frac{\alpha^{2}}{(\alpha+\lambda j)(\alpha+\lambda i)}-\sum_{j=1}^{n}\bigg(\frac{\alpha}{\alpha+\lambda j}\bigg)^{2}-2\sum_{i<j}\frac{\alpha^{2}}{(\alpha+\lambda j)(\alpha+\lambda i)}\nonumber\\
&=\ \sum_{j=1}^{n}\frac{\alpha}{\alpha+\lambda j}\bigg(1-\frac{\alpha}{\alpha+\lambda j}\bigg). \nonumber
\end{align}
\begin{lemma}
	For $n\ge 1$, let $\mu_{n,\lambda}(\alpha)$ be the mean of $Z_{n,\lambda}(\alpha)$, and let $\sigma_{n,\lambda}^{2}$ be the variance of $Z_{n,\lambda}(\alpha)$. Then we have 
	\begin{displaymath}
		\mu_{n,\lambda}(\alpha)=\sum_{j=1}^{n}\frac{\alpha}{\alpha+\lambda j},
	\end{displaymath}
	and 
	\begin{displaymath}
		\sigma_{n,\lambda}^{2}(\alpha)=\sum_{j=1}^{n}\frac{\alpha}{\alpha+\lambda j}\bigg(1-\frac{\alpha}{\alpha+\lambda j}\bigg). 
	\end{displaymath}
\end{lemma}
From \eqref{6}, we have 
\begin{equation}
\begin{aligned}
E\big[(1+\frac{z}{\lambda})^{X_{j}}\big]\ &=\ (1+\frac{z}{\lambda})P\{X_{j}=1\}+P\{X_{j}=0\}\\
&\ =\frac{1}{j}+\frac{z}{\lambda j}+1-\frac{1}{j}=1+\frac{z}{\lambda j}	.
\end{aligned}
\label{16}
\end{equation}
By \eqref{16}, we get 
\begin{equation}
\begin{aligned}
	E\Big[(1+\frac{z}{\lambda})^{Y_{n}}\Big]\ &= E\Big[(1+\frac{z}{\lambda})^{\sum_{j=1}^{n}X_{j}}\Big] \\
	&=\ \prod_{j=1}^{n}E\big[(1+\frac {z}{\lambda})^{X_{j}}\big]=\prod_{j=1}^{n}\bigg(1+\frac{z}{\lambda j}\bigg). 
\end{aligned}\label{17}	
\end{equation}
On the other hand, 
\begin{equation}
E\Big[z^{Z_{n,\lambda}(\alpha)}\Big]=E\Big[z^{\sum_{j=1}^{n}X_{j,\lambda}(\alpha)}\Big]=\prod_{j=1}^{n}E\Big[z^{X_{j,\lambda}(\alpha)}\Big]. \label{18}
\end{equation}
Note that 
\begin{align}
E\Big[z^{X_{j,\lambda}(\alpha)}\Big]\ &=\ zP\{X_{j,\lambda}(\alpha)=1\}+P\{X_{j,\lambda}(\alpha)=0\}\label{19}\\
&=\ \frac{\alpha z}{\alpha+\lambda j}+1-\frac{\alpha}{\alpha+\lambda j}	\nonumber \\
&=\ \frac{\alpha z+\lambda j}{\alpha+\lambda j},\quad (j=1,2,\dots,n). \nonumber
\end{align}
By \eqref{18} and \eqref{19}, we get 
\begin{equation}
E\Big[z^{Z_{n,\lambda}(\alpha)}\Big]=\prod_{j=1}^{n}\bigg(\frac{\alpha z+\lambda j}{\alpha+\lambda j}\bigg)\label{20}	.
\end{equation}
From \eqref{17} and \eqref{20}, we note that 
\begin{align}
&E\Big[(1+\frac{\alpha}{\lambda})^{Y_{n}}\Big] E\Big[z^{Z_{n,\lambda}(\alpha)}\Big]\label{21}\\
&\quad =\ \prod_{j=1}^{n}\bigg(1+\frac{\alpha}{\lambda j}\bigg)\prod_{j=1}^{n}\bigg(\frac{\alpha z+\lambda j}{\alpha+\lambda j}\bigg)\nonumber \\
&\quad=\ \prod_{j=1}^{n}\bigg(\frac{\alpha z+\lambda j}{\lambda j}\bigg)
=\ \prod_{j=1}^{n}\bigg(1+\frac{\alpha z}{\lambda j}\bigg).\nonumber
\end{align}
By \eqref{6}, we easily get 
\begin{equation}
E\bigg[\bigg(1+\frac{\alpha z}{\lambda}\bigg)^{Y_{n}}\bigg]=\prod_{j=1}^{n}E\bigg[\bigg(1+\frac{\alpha z}{\lambda}\bigg)^{X_{j}}\bigg]=\prod_{j=1}^{n}\bigg(1+\frac{\alpha z}{\lambda j}\bigg). \label{22}	
\end{equation}
Therefore, by \eqref{21} and \eqref{22}, we obtain the following theorem. 
\begin{theorem}
For $n\in\mathbb{N}$, we have 
\begin{displaymath}
E\Big[(1+\frac{\alpha}{\lambda})^{Y_{n}}\Big]E\Big[z^{Z_{n,\lambda}(\alpha)}\Big]=E\bigg[\bigg(1+\frac{\alpha z}{\lambda}\bigg)^{Y_{n}}\bigg]. 
\end{displaymath}
\end{theorem}

\section{Applications to sums of products of degenerate Stirling numbers of the second kind and Stirling numbers of the first kind}
Here, as an application of the dimorphic property in Theorem 2, we derive two different expressions on the sum $\sum_{m=l}^{n}S_{2,\lambda}(n+1,m+1)S_{1}(m+1,l+1)$, one involving $Y_{n}$ (see Theorem 3) and the other involving $Z_{n,\lambda}(\alpha)$ (see Theorem 4). \par
From \eqref{2} and \eqref{4}, we note that 
\begin{align}
(x)_{n+1,\lambda}\ &=\ \sum_{m=0}^{n+1}S_{2,\lambda}(n+1,m)(x)_{m}\ =\ \sum_{m=1}^{n+1}S_{2,\lambda}(n+1,m)(x)_{m} \label{23} \\
&=\ \sum_{m=0}^{n}S_{2,\lambda}(n+1,m+1)(x)_{m+1} \nonumber\\
&=\ x\sum_{m=0}^{n}S_{2,\lambda}(n+1,m+1)\sum_{l=0}^{m}S_{1}(m+1,l+1)x^{l} \nonumber \\
&=\ x\sum_{l=0}^{n}\bigg(\sum_{m=l}^{n}S_{2,\lambda}(n+1,m+1)S_{1}(m+1,l+1)\bigg)x^{l}. \nonumber
\end{align}
By \eqref{1} and \eqref{17}, we easily get 
\begin{align}
(x)_{n+1,\lambda}\ &=\ x(x-\lambda)(x-2\lambda)\cdots (x-n\lambda)\label{24} \\
&=\ n!\lambda^{n}(-1)^{n}x\bigg(1-\frac{x}{\lambda}\bigg)\bigg(1-\frac{x}{\lambda 2}\bigg)\cdots\bigg(1-\frac{x}{\lambda n}\bigg)\nonumber \\
&=\ n!(-1)^{n}\lambda^{n}x\prod_{j=1}	^{n}\bigg(1-\frac{x}{\lambda j}\bigg)\nonumber\\
&=\ n!(-1)^{n}\lambda^{n}xE\bigg[\bigg(1-\frac{x}{\lambda})^{Y_{n}}\bigg].\nonumber
\end{align}
From \eqref{23} and \eqref{24}, we note that 
\begin{align}
&z\sum_{l=0}^{n}\bigg(\sum_{m=l}^{n}S_{2,\lambda}(n+1,m+1)S_{1}(m+1,l+1)\bigg)z^{l}\label{25}\\
&\quad=n!(-1)^{n}\lambda^{n}zE\bigg[\bigg(1-\frac{z}{\lambda}\bigg)^{Y_{n}}\bigg]	\nonumber \\
&\quad=n!(-1)^{n}\lambda^{n}z\sum_{l=0}^{n}(-1)^{l}E\bigg[\binom{Y_{n}}{l}\bigg]\bigg(\frac{z}{\lambda}\bigg)^{l}\nonumber \\
&\quad=z\sum_{l=0}^{n}(-1)^{n-l}\lambda^{n-l}n!E\bigg[\binom{Y_{n}}{l}\bigg]z^{l}.\nonumber
\end{align}
By comparing the coefficients on both sides of \eqref{25}, we obtain the following theorem. 
\begin{theorem}
For $0\le l\le n$, we have 
\begin{displaymath}
\sum_{m=l}^{n}S_{2,\lambda}(n+1,m+1)S_{1}(m+1,l+1)=(-1)^{n-l}\lambda^{n-l}n!E\bigg[\binom{Y_{n}}{l}\bigg].
\end{displaymath}
\end{theorem}
As $Y_{n}$ is taking integer values between $0$ and $n$, we see that 
\begin{align}
\frac{1}{2\pi}\int_{-\pi}^{\pi}E\Bigg[\bigg(1+\frac{\alpha}{\lambda}e^{i\theta}\bigg)^{Y_{n}}\Bigg]e^{-im\theta}d\theta\ &=\ \frac{1}{2\pi} \int_{-\pi}^{\pi}
\sum_{j=0}^{n}E\bigg[\binom{Y_{n}}{j}\bigg]\bigg(\frac{\alpha	}{\lambda}\bigg)^{j}e^{ji\theta}e^{-im\theta}d\theta \label{26} \\
&=\ \sum_{j=0}^{n}E\bigg[\binom{Y_{n}}{j}\bigg]\frac{\alpha^{j}}{2\pi}\lambda^{-j} \int_{-\pi}^{\pi}e^{i\theta(j-m)}d\theta \nonumber \\
&=\lambda^{-m}\alpha^{m}E\bigg[\binom{Y_{n}}{m}\bigg],\nonumber
\end{align}
where $i=\sqrt{-1}$. \par 
From Theorem 2 and \eqref{26}, we have 
\begin{align}
&\lambda^{-l}\alpha^{l}E\bigg[\binom{Y_{n}}{l}\bigg]=\frac{1}{2\pi} \int_{-\pi}^{\pi}E\bigg[\bigg(1+\frac{\alpha}{\lambda}e^{i\theta}\bigg)^{Y_{n}}\bigg]e^{-il\theta}d\theta	\label{27} \\
&\quad=\ \frac{1}{2\pi}E\Big[(1+\frac {\alpha}{\lambda})^{Y_{n}}\Big] \int_{-\pi}^{\pi}E\Big[e^{i\theta Z_{n,\lambda}(\alpha)}\Big]e^{-il\theta}d\theta \nonumber \\
&\quad=\ \prod_{j=1}^{n}\bigg(1+\frac{\alpha}{\lambda j}\bigg)\frac{1}{2\pi} \int_{-\pi}^{\pi}E\Big[e^{i\theta(Z_{n,\lambda}(\alpha)-l)}\Big]d\theta \nonumber \\
&\quad=\ \frac{\Gamma(1+\frac{\alpha}{\lambda})(1+\frac{\alpha}{\lambda})(2+\frac{\alpha}{\lambda})\cdots(n+\frac{\alpha}{\lambda})}{n!\Gamma(1+\frac{\alpha}{\lambda})}\frac{1}{2\pi} \int_{-\pi}^{\pi}E\Big[e^{i\theta(Z_{n,\lambda}(\alpha)-l)}\Big]d\theta\nonumber \\
&\quad=\frac{\Gamma(n+\frac{\alpha}{\lambda}+1)}{n!\Gamma(1+\frac{\alpha}{\lambda})}\frac{1}{2\pi} \int_{-\pi}^{\pi} E\Big[e^{i\theta(Z_{n,\lambda}(\alpha)-l)}\Big]d\theta\nonumber .
\end{align}
By Theorem 3 and \eqref{27}, we get 
\begin{align}
&\sum_{m=l}^{n}S_{2,\lambda}(n+1,m+1)S_{1}(m+1,l+1)\ =\ (-1)^{n-l}\lambda^{n-l}n!E\bigg[\binom{Y_{n}}{l}\bigg]\label{28} \\
&\quad=\ (-1)^{n-l}\lambda^{n}\frac{\Gamma(n+\frac{\alpha}{\lambda}+1)}{\alpha^{l}\Gamma(1+\frac{\alpha}{\lambda})} \frac{1}{2\pi} \int_{-\pi}^{\pi} \Big[e^{i\theta(Z_{n,\lambda}(\alpha)-l)}\Big]d\theta\nonumber \\
&\quad=\ (-1)^{n-l}\lambda^{n+1}\frac{\Gamma(n+\frac{\alpha}{\lambda}+1)}{\alpha^{l+1}\Gamma(\frac{\alpha}{\lambda})} \frac{1}{2\pi} \int_{-\pi}^{\pi} \Big[e^{i\theta(Z_{n,\lambda}(\alpha)-l)}\Big]d\theta\nonumber.
\end{align}
Therefore, by \eqref{28}, we obtain the following theorem. 
\begin{theorem}
For $0\le l\le n$, $\alpha>0 (\alpha\in\mathbb{R})$, $\lambda\in (0,1)$, and $i=\sqrt{-1}$, we have 
\begin{displaymath}
\frac{1}{2\pi} \int_{-\pi}^{\pi} \Big[e^{i\theta(Z_{n,\lambda}(\alpha)-l)}\Big]d\theta=(-1)^{n-l}\frac{\alpha^{l+1}}{\lambda^{n+1}}\frac{\Gamma(\frac{\alpha}{\lambda})}{\Gamma(n+\frac{\alpha}{\lambda}+1)}\sum_{m=l}^{n}S_{2,\lambda}(n+1,m+1)S_{1}(m+1,l+1). 
\end{displaymath}
\end{theorem}
\section{Conclusion}
Let $(X_{j})_{1\le j\le n}$ be identically independent Bernoulli random variables such that $X_{j}$ has the probability of success $\frac{1}{j},\ (j=1,2,\dots,n)$, with $Y_{n}=X_{1}+X_{2}+\cdots+X_{n}$. 	
Let $\big(X_{j,\lambda}(\alpha)\big)_{1\le j\le n}$ be identically independent Bernoulli random variables such that $X_{j,\lambda}(\alpha)$ has the probability of success $\frac{\alpha}{\alpha+\lambda j},\ (j=1,2,\dots,n)$, with $Z_{n,\lambda}(\alpha)=X_{1,\lambda}(\alpha)+X_{2,\lambda}(\alpha)+\cdots+X_{n,\lambda}(\alpha)$. Here $\alpha$ is a positive real number and $\lambda$ is a real number with $0<\lambda<1$. \par
Then we showed in Theorem 2 the dimorphic property:
\begin{displaymath}
E\Big[(1+\frac{\alpha}{\lambda})^{Y_{n}}\Big]E\Big[z^{Z_{n,\lambda}(\alpha)}\Big]=E\bigg[\bigg(1+\frac{\alpha z}{\lambda}\bigg)^{Y_{n}}\bigg]. 
\end{displaymath}
Further, we derived two different expressions on the sum $\sum_{m=l}^{n}S_{2,\lambda}(n+1,m+1)S_{1}(m+1,l+1)$ in connection with $Y_{n}$ and $Z_{n,\lambda}(\alpha)$. Indeed, we derived the following:
\begin{align*}
&\sum_{m=l}^{n}S_{2,\lambda}(n+1,m+1)S_{1}(m+1,l+1)\ =\ (-1)^{n-l}\lambda^{n-l}n!E\bigg[\binom{Y_{n}}{l}\bigg]\\
&\quad=\ (-1)^{n-l}\lambda^{n+1}\frac{\Gamma(n+\frac{\alpha}{\lambda}+1)}{\alpha^{l+1}\Gamma(\frac{\alpha}{\lambda})} \frac{1}{2\pi} \int_{-\pi}^{\pi} \Big[e^{i\theta(Z_{n,\lambda}(\alpha)-l)}\Big]d\theta\nonumber.
\end{align*} \par
There are various ways of studying special numbers and polynomials, to mention a few, generating functions, combinatorial methods, probability theory, $p$-adic analysis, umbral calculus, differential equations, special functions and analytic number theory. 
In recent years, we have had lively interests in the study of various degenerate versions of special numbers and polynomials with those diverse tools.
As a result of such explorations, we came up with, for example, the degenerate Stirling numbers which are degenerate versions of the ordinary Stirling numbers and appear in many different contexts.
The novelty of this paper is that we obtained two different expressions on sums of products of degenerate Stirling numbers of the second kind and Stirling numbers of the first kind in connection with two different sums of identically independent Bernoulli random variables.
This is one example of our efforts in the applications of probability theory to the study of some special numbers and polynomials and also of degenerate versions of those numbers and polynomials. \\
\indent We would like to continue to find many applications of probability theory to the study of some special numbers and polynomials and also of their degenerate versions. More generally, it is one of our future projects to continue to explore various degenerate versions of many special polynomials and numbers with aforementioned tools.

\end{document}